\documentclass{article}
\usepackage{amsmath, amssymb, graphics}

\usepackage{amsthm}
\usepackage[dvipdfmx]{graphicx}

\newtheorem{thm}{Theorem}[section]
\newtheorem{df}{Definition}[section]
\newtheorem{rem}{Remark}[section]

\newtheorem{prp}{Proposition}[section]
\newtheorem{lem}{Lemma}[section]
\newtheorem{cor}{Corollary}[section]

\newcommand{\mathsym}[1]{{}}
\newcommand{\unicode}[1]{{}}

\makeatletter
    
    \@addtoreset{equation}{section}
  \makeatother
  
\begin{document}
\title{Iterative method of construction of good rhythms}
\author{Fumio HAZAMA\\Tokyo Denki University\\Hatoyama, Hiki-Gun, Saitama JAPAN\\
e-mail address:hazama@mail.dendai.ac.jp\\Phone number: (81)49-296-5267}
\date{\today}
\maketitle
\thispagestyle{empty}

\begin{abstract}
On the space of rhythms of arbitrary length with a fixed number of onsets, a self map $F$ is constructed. It is shown that for any rhythm $\mathbf{r}$ of the space there exists a nonnegative integer $k$ such that $F^k(\mathbf{r})$ falls into the category of "maximally even rhythms" investigated in [1]. 
\end{abstract}

\section{Introduction}
In [1], the authors introduce the reader to several methods of construction of "good" rhythms, and show how various mathematical tools play crucial roles for their investigation. In this paper we propose a method which enables us to construct good rhythms {\it iteratively}. For any integer $N\geq 3$ and for any integer $n<N$, we employ the set ${\rm CR}_N^{(n)}=\{\mathbf{b}\in\mathbf{Z}_2^N;d(\mathbf{b},\mathbf{0})=n\}$ as a representation space of the set of rhythms with period $N$ and with $n$ onsets, where $d(\cdot,\cdot)$ denote the Hamming distance. The main purpose of this paper is to construct a self-map $F:{\rm CR}_N^{(n)}\rightarrow {\rm CR}_N^{(n)}$ such that, for any given rhythm $\mathbf{r}\in {\rm CR}_N^{(n)}$, there exists a nonnegative integer $K$ with the property that $F^K(\mathbf{r})$ becomes {\it good}. What we mean by a good rhythm will be explained in the body of this paper. At this point we merely say that our good rhythms contain those exhibited in [1] and some other interesting ones.

\section{Cyclic polygons modulo $N$}
As is described in [1], a rhythm can be specified in various ways. The set ${\rm CR}_N^{(n)}$ is basic in the sense that it specifies directly when we sound $n$ onsets among $N$ pulses by the places of nonzero entries. We, however, attach a cyclic polygon to a rhythm which will allow us to investigate the iterative process, mentioned in Introduction, more visually. \\

Let $\mu_N$ denote the set of the $N$-th roots of unity in the complex plane $\mathbf{C}$, and let $\chi_N:\mathbf{Z}\rightarrow\mu_N$ denote the map defined by $\chi_N(k)=\zeta_N^k$ with $\zeta_N=\exp(2\pi i/N)\in\mu_N$. In this paper we denote the set $[0,N-1]=\{0,1,\cdots,N-1\}$ by $\mathbf{Z}_N$ and provide it with a group operation "$+_N$" defined by the rule
\begin{eqnarray*}
a+_Nb=R_N(a+b)
\end{eqnarray*}
where $R_N:\mathbf{Z}\rightarrow\mathbf{Z}_N$ is the map which takes the remainder of an integer divided by $N$. The reason why we fix a set of representatives is the fact that we employ the natural total order on $[0,N-1]$ inherited from $\mathbf{Z}$ later in this paper. For a positive integer $n\geq 3$, an $n$-gon $\mathbf{P}=(P_0,P_1,\cdots, P_{n-1})$ is said to be a {\it cyclic} $n$-{\it gon mod} $N$ if every vertex belongs to $\mu_N$ and the vertices are numbered counterclockwise. We assume that $P_i\neq P_{i+1}$ for every $i\in [0, n-1]$ and $P_n$ is understood to be equal to $P_0$. We denote the set of cyclic $n$-gons mod $N$ by ${\rm CP}_N^{(n)}$. Note that each vertex $P_i, i\in[0,n-1],$ can be expressed uniquely as $\chi_N(a_i)$ with $a_i\in\mathbf{Z}_N$. Then the sequence $\mathbf{a}=(a_0,\cdots,a_{n-1})$ has the following property (A):\\

\noindent
(A): There exists a unique number $k_0\in[0,n-1]$ such that $a_{k_0}>a_{k_0+1}$ and that $a_k<a_{k+1}$ holds for every other number $k\in[0,n-1]\setminus\{k_0\}$.\\

\noindent
We call the number $k_0$ stated in (A) the {\it jumping number} of the cyclic polygon $\mathbf{P}$ and denote it by $jump(\mathbf{P})$. Furthermore we denote the set of sequences $\mathbf{a}=(a_0,\cdots,a_{n-1})\in\mathbf{Z}_N^n$ which satisfy the property (A) by ${\rm CA}_N^{(n)}$, and define a map ${\rm X}_N:{\rm CA}_N^{(n)}\rightarrow {\rm CP}_N^{(n)}$ by the rule
\begin{eqnarray*}
{\rm X}_N(a_0,\cdots,a_{n-1})=(\chi_N(a_0),\cdots, \chi_N(a_{n-1})).
\end{eqnarray*}
Noticing that ${\rm X}_N$ is a bijection, we call $jump({\rm X}_N(\mathbf{a}))$ also the {\it jumping number} of $\mathbf{a}$ and denote it simply by $jump(\mathbf{a})$.

\section{Discrete averages}
In this section we introduce two types of discrete averages, one for the elements of $\mu_N$, and one for those of $\mathbf{Z}_N$, and investigate several properties which they share.

\begin{df}
For any $P=\chi_N(a),a\in\mathbf{Z}_N,$ let $sq_N:\mu_N\rightarrow \mu_N$ denote the map defined by 
\begin{eqnarray*}
sq_N(P)=\chi_N(\left\lfloor\frac{a}{2}\right\rfloor),
\end{eqnarray*}
and call it the {\rm discrete square root} of $P$. Furthermore for any pair $(P,Q)\in\mu_N\times\mu_N$, let
\begin{eqnarray*}
av_{\mu_N}(P,Q)=P\cdot sq_N(P^{-1}\cdot Q),
\end{eqnarray*}
where the dot "$\cdot$" means the usual product of the complex numbers, and call it the $\mu${\rm -discrete average}.
\end{df}

\noindent
Note that the $\mu$-discrete average is sensitive to the order of $P,Q$, namely we have $av_{\mu_N}(P,Q)\neq av_{\mu_N}(Q,P)$ in general. Another type of discrete average is defined as follows:

\begin{df}
For any pair $(a,b)\in\mathbf{Z}_N\times\mathbf{Z}_N$, let
\begin{eqnarray*}
av_{\mathbf{Z}_N}(a,b)=a+_N\left\lfloor\frac{-a+_Nb}{2}\right\rfloor,
\end{eqnarray*}
and call it the $\mathbf{Z}${\rm -discrete average}.
\end{df}

In order to describe where the average falls into, it is convenient to introduce the following notation:

\begin{df}
For any $a,b\in\mathbf{Z}_N$, let
\begin{eqnarray*}
[a,b]=\left\{
\begin{array}{ll}
\{a,a+1,\cdots,b\}, & \mbox{if $a\leq b$},\\
\{a,a+1,\cdots,0,1,\cdots,b\}, & \mbox{if $a>b$}.
\end{array}
\right.
\end{eqnarray*}
Furthermore we define $[a,b)$ to be the set $[a,b]\setminus\{b\}$, as is in the case for the usual half-closed interval on the real line.
\end{df}

\noindent
Accordingly, we define the notion of an interval in $\mu_N$ as follows:

\begin{df}
For any $P,Q\in\mu_N$, let $P=\chi_N(a), Q=\chi_N(b)$ with $a,b\in\mathbf{Z}_N$. We define the sets $[P,Q], [P,Q)\subset \mu_N$ by the following rules respectively:
\begin{eqnarray*}
&&[P,Q]=\{\chi_N(x);x\in [a,b]\},\\
&&[P,Q)=\{\chi_N(x);x\in [a,b)\}.
\end{eqnarray*}
\end{df}

\begin{rem}
The interval $[P,Q]$ is also expressed as the set of points in $\mu_N$ which begins at $P$ and goes around the unit circle counterclockwise ending at $Q$. 
\end{rem}

\noindent
It follows from the definition that
\begin{eqnarray}
av_{\mathbf{Z}_N}(a,b)\in [a,b) 
\end{eqnarray}
holds for any $a,b\in\mathbf{Z}_N$, and that
\begin{eqnarray}
av_{\mu_N}(P,Q)\in [P,Q) 
\end{eqnarray}
holds for any $P,Q\in\mu_N$. As is expected, these two averages are related naturally as follows:

\begin{prp}
For any $a,b\in\mathbf{Z}_N$, we have
\begin{eqnarray}
\chi_N(av_{\mathbf{Z}_N}(a,b))=av_{\mu_N}(\chi_N(a),\chi_N(b)).
\end{eqnarray}
\end{prp}

\noindent
{\it Proof}. The left hand side of (2.3) is equal to
\begin{eqnarray*}
\chi_N(av_{\mathbf{Z}_N}(a,b))=\chi_N(a+_N\left\lfloor\frac{-a+_Nb}{2}\right\rfloor)
\end{eqnarray*}
by the very definition. The right hand side of (2.3) is computed as
\begin{eqnarray*}
av_{\mu_N}(\chi_N(a),\chi_N(b))&=&\chi_N(a)\cdot sq_N(\chi_N(a)^{-1}\cdot\chi_N(b))\\
&=&\chi_N(a)\cdot sq_N(\chi_N(-a+_Nb))\\
&=&\chi_N(a)\cdot \chi_N(\left\lfloor\frac{-a+_Nb}{2}\right\rfloor)\\
&=&\chi_N(a+_N\left\lfloor\frac{-a+_Nb}{2}\right\rfloor),
\end{eqnarray*}
hence the both sides are the same. \qed\\

\noindent
Based on the two averages, we introduce the following two types of average transformations:

\begin{df}
$(1)$ For any $\mathbf{a}=(a_0,a_1,\cdots,a_{n-1})\in{\rm CA}_N^{(n)}$, we define
\begin{eqnarray*}
{\rm dav}_A(\mathbf{a})=(av_{\mathbf{Z}_N}(a_0,a_1),av_{\mathbf{Z}_N}(a_1,a_2),\cdots,av_{\mathbf{Z}_N}(a_{n-1},a_0)),
\end{eqnarray*}
and call the map ${\rm dav}_A:{\rm CA}_N^{(n)}\rightarrow{\rm CA}_N^{(n)}$ the $A${\rm -discrete average transformation}.\\
$(2)$ For any $\mathbf{P}=(P_0,P_1,\cdots,P_{n-1})\in{\rm CP}_N^{(n)}$, we define
\begin{eqnarray*}
{\rm dav}_P(\mathbf{P})=(av_{\mu_N}(P_0,P_1),av_{\mu_N}(P_1,P_2),\cdots,av_{\mu_N}(P_{n-1},P_0)),
\end{eqnarray*}
and call the map ${\rm dav}_P:{\rm CP}_N^{(n)}\rightarrow{\rm CP}_N^{(n)}$ the $P${\rm -discrete average transformation}.
\end{df}

\begin{rem}
It follows from $(2.1)$ and $(2.2)$ that the image of the map ${\rm dav}_P$ (resp. ${\rm dav}_A$) is contained in ${\rm CP}_N^{(n)}$ (resp. ${\rm CA}_N^{(n)}$) so that our use of the term "average {\it transformation}" is justified.
\end{rem}

As is imagined from Proposition 2.1, the two average transformations are related through the character map ${\rm X}_N$:

\begin{prp}
For any $\mathbf{a}\in{\rm CA}_N^{(n)}$, we have the equality
\begin{eqnarray*}
{\rm X}_N({\rm dav}_A(\mathbf{a}))={\rm dav}_P({\rm X}_N(\mathbf{a})).
\end{eqnarray*}
\end{prp}

\noindent
{\it Proof}. Let $\mathbf{a}=(a_0,a_1,\cdots, a_{n-1})$. Then we can compute as follows:
\begin{eqnarray*}
{\rm X}_N({\rm dav}_A(\mathbf{a}))&=&{\rm X}_N({\rm dav}_A(a_0,a_1,\cdots,a_{n-1}))\\
&=&{\rm X}_N(av_{\mathbf{Z}_N}(a_0,a_1), \cdots, av_{\mathbf{Z}_N}(a_{n-1},a_0))\\
&=&(\chi_N(av_{\mathbf{Z}_N}(a_0,a_1)), \cdots, \chi_N(av_{\mathbf{Z}_N}(a_{n-1},a_0)))\\
&=&(av_{\mu_N}(\chi_N(a_0),\chi_N(a_1)), \cdots, av_{\mu_N}(\chi_N(a_{n-1}),\chi_N(a_0)))\\
&&\hspace{40mm}(\Leftarrow \mbox{by Proposition 2.1})\\
&=&{\rm dav}_P(\chi_N(a_0), \cdots, \chi_N(a_{n-1}))\\
&=&{\rm dav}_P({\rm X}_N(\mathbf{a})).
\end{eqnarray*}
This completes the proof. \qed\\

\section{Sequences of differences}
When we try to understand the behavior of the iteration of the average transformation ${\rm dav}_A$, it is helpful for us to look at the sequence of differences of an element of ${\rm CA}_N^{(n)}$.

\begin{df}
For any $\mathbf{a}=(a_0,\cdots,a_{n-1})\in{\rm CA}_N^{(n)}$, let
\begin{eqnarray*}
{\rm diff}(\mathbf{a})=(b_0,\cdots,b_{n-1}),
\end{eqnarray*}
where
\begin{eqnarray*}
b_k=a_{k+1}-_Na_k,\hspace{2mm}k\in [0,n-1],
\end{eqnarray*}
and call the value ${\rm diff}(\mathbf{a})$ the A-${\rm difference}$ of $\mathbf{a}$.
\end{df} 

In order to specify the range of the $A$-difference, we introduce the following notation:
\begin{eqnarray}
{\rm CD}_N^{(n)}=\{(b_0,\cdots,b_{n-1})\in\mathbf{Z}_N^n;\sum_{k=0}^{n-1}b_k=N\}.
\end{eqnarray}
Note here that the summation appearing on the right hand side of (3.1) should be understood as that of integers in $\mathbf{Z}$. 

\begin{prp}
For any $\mathbf{a}=(a_0,\cdots,a_{n-1})\in{\rm CA}_N^{(n)}$, we have ${\rm diff}(\mathbf{a})\in{\rm CD}_N^{(n)}$. Therefore the A-difference defines actually a map ${\rm diff}:{\rm CA}_N^{(n)}\rightarrow {\rm CD}_N^{(n)}$.
\end{prp}

\noindent
{\it Proof}. Let $jump(\mathbf{a})=k_0$ and let ${\rm diff}(\mathbf{a})=\mathbf{b}=(b_0,\cdots,b_{n-1})$. Then we have
\begin{eqnarray*}
b_{k_0}=a_{k_0+1}-_Na_{k_0}=a_{k_0+1}-a_{k_0}+N.
\end{eqnarray*}
On the other hand, for any number $k\in [0,n-1]\setminus\{k_0\}$ we have $0<a_{k+1}-a_k\leq N-1$, and hence
\begin{eqnarray*}
a_{k+1}-_Na_k=a_{k+1}-a_k.
\end{eqnarray*}
Therefore we can compute as follows:
\begin{eqnarray*}
\sum_{k=0}^{n-1}b_k&=&\sum_{k=0}^{k_0-1}b_k+b_{k_0}+\sum_{k=k_0+1}^{n-1}b_k\\
&=&\sum_{k=0}^{k_0-1}(a_{k+1}-a_k)+(a_{k_0+1}+N-a_{k_0})+\sum_{k=k_0+1}^{n-1}(a_{k+1}-a_k)\\
&=&(a_{k_0}-a_0)+(a_{k_0+1}+N-a_{k_0})+(a_{n}-a_{k_0+1})\\
&=&N+a_n-a_0\\
&=&N\hspace{5mm}(\Leftarrow \mbox{since $a_n=a_0$})
\end{eqnarray*}
This completes the proof. \qed\\

\noindent

We examine a few of examples in order to observe the interdependence of the $A$-average and the $A$-difference. For ease of description we put
\begin{eqnarray*}
\mathbf{a}^{(k)}&=&{\rm dav}_A^k(\mathbf{a}),\\
\mathbf{d}^{(k)}&=&{\rm diff}(\mathbf{a}^{(k)}), \\
\mathbf{P}^{(k)}&=&{\rm X}_N(\mathbf{a}^{(k)})
\end{eqnarray*}
for any nonnegative integer $k$.\\

\noindent
Example 3.1. When $N=16, n=4$ (so that $n$ divides $N$), and $\mathbf{a}=(0,3,7,14)$, we can compute them as follows:
\begin{eqnarray*}
\begin{array}{lcl}
\mathbf{a}^{(0)}=(0,3,7,14) & \rightarrow & \mathbf{d}^{(0)}=(3,4,7,2) \\
\downarrow & & \\
\mathbf{a}^{(1)}=(1,5,10,15) & \rightarrow & \mathbf{d}^{(1)}=(4,5,5,2) \\
\downarrow & & \\
\mathbf{a}^{(2)}=(3,7,12,0) & \rightarrow & \mathbf{d}^{(2)}=(4,5,4,3) \\
\downarrow & & \\
\mathbf{a}^{(3)}=(5,9,14,1) & \rightarrow & \mathbf{d}^{(3)}=(4,5,3,4) \\
\downarrow & & \\
\mathbf{a}^{(4)}=(7,11,15,3) & \rightarrow & \mathbf{d}^{(4)}=(4,4,4,4) \\
\downarrow & & \\
\mathbf{a}^{(5)}=(9,13,1,5) & \rightarrow & \mathbf{d}^{(5)}=(4,4,4,4) \\
\downarrow & & \\
\mathbf{a}^{(6)}=(11,15,3,7) & \rightarrow & \mathbf{d}^{(6)}=(4,4,4,4) \\
\downarrow & & \\
\mathbf{a}^{(7)}=(13,1,5,9) & \rightarrow & \mathbf{d}^{(7)}=(4,4,4,4) \\
\vdots &  & \hspace{10mm}\vdots \\
\end{array}
\end{eqnarray*}
As is seen from this, the differences $\mathbf{d}^{(k)}$ for $k\geq 4$ stabilize to $(4,4,4,4)$, therefore  every polygon $\mathbf{P}^{(k)}$ becomes the square for $k\geq 4$.\\

\noindent
Example 3.2. When $N=16, n=5$, and $\mathbf{a}=(0,1,3,7,14)$, the successive averages and differences become as follows:
\begin{eqnarray*}
\begin{array}{lcl}
\mathbf{a}_0=(0,1,3,7,14) & \rightarrow & \mathbf{d}_0=(1,2,4,7,2) \\
\downarrow & & \\
\mathbf{a}_1=(0,2,5,10,15) & \rightarrow & \mathbf{d}_1=(2,3,5,5,1) \\
\downarrow & & \\
\mathbf{a}_2=(1,3,7,12,15) & \rightarrow & \mathbf{d}_2=(2,4,5,3,2) \\
\downarrow & & \\
\mathbf{a}_3=(2,5,9,13,0) & \rightarrow & \mathbf{d}_3=(3,4,4,3,2) \\
\downarrow & & \\
\mathbf{a}_4=(3,7,11,14,1) & \rightarrow & \mathbf{d}_4=(4,4,3,3,2) \\
\downarrow & & \\
\mathbf{a}_5=(5,9,12,15,2) & \rightarrow & \mathbf{d}_5=(4,3,3,3,3) \\
\downarrow & & \\
\mathbf{a}_6=(7,10,13,0,3) & \rightarrow & \mathbf{d}_6=(3,3,3,3,4) \\
\downarrow & & \\
\mathbf{a}_7=(8,11,14,1,5) & \rightarrow & \mathbf{d}_7=(3,3,3,4,3) \\
\downarrow & & \\
\mathbf{a}_8=(9,12,15,3,6) & \rightarrow & \mathbf{d}_8=(3,3,4,3,3) \\
\downarrow & & \\
\mathbf{a}_9=(10,13,1,4,7) & \rightarrow & \mathbf{d}_9=(3,4,3,3,3) \\
\downarrow & & \\
\mathbf{a}_{10}=(11,15,2,5,8) & \rightarrow & \mathbf{d}_{10}=(4,3,3,3,3) \\
\vdots &  & \hspace{10mm}\vdots \\
\end{array}
\end{eqnarray*}
In this case $N$ and $n$ are coprime so that the corresponding polygons $\mathbf{P}^{(k)}$ cannot be regular pentagons for any $k\geq 0$. The differences $\mathbf{d}^{(k)}$ for $k\geq 5$, however, exhibit a periodic behavior.\\

The main purpose of this paper is to show that the phenomena in these two examples are typical and can be generalized in a natural way.

\section{$D$-average}
In the last two examples we compute $\mathbf{d}^{(k+1)}$ via the average transformation of $\mathbf{a}^{(k)}$ for each $k$. In this section we show that we can compute the next $\mathbf{d}^{(k+1)}$ from $\mathbf{d}^{(k)}$ directly.\\

\subsection{Definition and fundamental properties of $D$-average}
The following function plays a crucial role for the definition of $D$-average:

\begin{df}
For any $a,b\in\mathbf{Z}$, we define 
\begin{eqnarray*}
av_{fc}(a,b)=
\left\{
\begin{array}{ll}
\left\lfloor\frac{a+b}{2}\right\rfloor, & \mbox{if $a\equiv 0\pmod 2$},\\
\left\lceil\frac{a+b}{2}\right\rceil, & \mbox{if $a\equiv 1\pmod 2$},\\
\end{array}
\right.
\end{eqnarray*}
and call the resulting map $av_{fc}:\mathbf{Z}^2\rightarrow\mathbf{Z}$ the fc{\rm-average} of $a$ and $b$.
\end{df}

\begin{rem}
The name $fc$-{\it average} comes from the fact that it employs the {\underline f}loor or {\underline c}eiling function in this order according as the first argument $a$ is even or odd.
\end{rem}

\noindent
The $D$-average is introduced as follows:

\begin{df}
For any sequence $\mathbf{d}=(d_0,\cdots, d_{n-1})\in\mathbf{Z}^n$, we put
\begin{eqnarray*}
{\rm dav}_{fc}(\mathbf{d})=(av_{fc}(d_0,d_1),av_{fc}(d_1,d_2),\cdots,av_{fc}(d_{n-1},d_0)),
\end{eqnarray*}
and call the map ${\rm dav}_{fc}:\mathbf{Z}^n\rightarrow\mathbf{Z}^n$ the $D${\rm -discrete average transformation}.\\
\end{df}

The following proposition shows that the $D$-average ${\rm dav}_{fc}$ maps ${\rm CD}_N^{(n)}$ into itself. In order to show this fact, we put 
\begin{eqnarray*}
s(\mathbf{d})=\sum_{i=0}^{n-1}d_i,
\end{eqnarray*}
for any $\mathbf{d}=(d_0,\cdots,d_{n-1})\in \mathbf{Q}^n$. Then we have the following:

\begin{prp}
For any $\mathbf{d}\in\mathbf{Z}^n$, the equality
\begin{eqnarray*}
s({\rm dav}_{fc}(\mathbf{d}))=s(\mathbf{d})
\end{eqnarray*}
holds. In particular, the $D$-average transformation defines a map
\begin{eqnarray*}
{\rm dav}_{fc}:{\rm CD}_N^{(n)}\rightarrow {\rm CD}_N^{(n)}.
\end{eqnarray*}
\end{prp}

\noindent
{\it Proof}. Let $av:\mathbf{Q}^2\rightarrow\mathbf{Q}$ denote the usual average defined by $av(a,b)=(a+b)/2$, and let ${\rm av}:\mathbf{Q}^n\rightarrow\mathbf{Q}^n$ denote the map defined by 
\begin{eqnarray*}
{\rm av}(\mathbf{d})=(av(d_0,d_1),\cdots,av(d_{n-1},d_0))
\end{eqnarray*}
for any $\mathbf{d}=(d_0,\cdots, d_{n-1})\in\mathbf{Q}^n$. Then we have
\begin{eqnarray}
s({\rm av}(\mathbf{d}))&=&\frac{d_0+d_1}{2}+\cdots+\frac{d_{n-1}+d_0}{2}\nonumber\\
&=&\frac{1}{2}(2d_0+\cdots 2d_{n-1})\nonumber\\
&=&d_0+\cdots +d_{n-1}\nonumber\\
&=&s(\mathbf{d}).
\end{eqnarray}
On the other hand, we have
\begin{eqnarray*}
av_{fc}(d_k,d_{k+1})=
\left\{
\begin{array}{ll}
av(d_k,d_{k+1}), & \mbox{if }(d_k,d_{k+1})\equiv (0,0)\pmod 2,\\
av(d_k,d_{k+1}), & \mbox{if }(d_k,d_{k+1})\equiv (1,1)\pmod 2,\\
av(d_k,d_{k+1})-\frac{1}{2}, & \mbox{if} (d_k,d_{k+1})\equiv (0,1)\pmod 2,\\
av(d_k,d_{k+1})+\frac{1}{2}, & \mbox{if} (d_k,d_{k+1})\equiv (1,0)\pmod 2.
\end{array}
\right.
\end{eqnarray*}
Hence if we put
\begin{eqnarray*}
n_{(0,1)}&=&|\{(d_k,d_{k+1});k\in [0,n-1], (d_k,d_{k+1})\equiv (0,1)\pmod 2\}|,\\
n_{(1,0)}&=&|\{(d_k,d_{k+1});k\in [0,n-1], (d_k,d_{k+1})\equiv (1,0)\pmod 2\}|.
\end{eqnarray*}
then we see that the equality
\begin{eqnarray}
\sum_{k=0}^{n-1}av_{fc}(d_k,d_{k+1})=\sum_{k=0}^{n-1}av(d_k,d_{k+1})-\frac{1}{2}n_{(0,1)}+\frac{1}{2}n_{(1,0)}
\end{eqnarray}
holds. Since the sequence $(d_0,\cdots,d_{n-1},d_0)$ either begins with a even number and ends with a even number, or begins with an odd number and ends with an odd number, we must have $n_{(0,1)}=n_{(1,0)}$. Therefore the second and the third terms on the right hand side of (4.2) cancels each other, and we have
\begin{eqnarray*}
\sum_{k=0}^{n-1}av_{fc}(d_k,d_{k+1})=\sum_{k=0}^{n-1}av(d_k,d_{k+1}).
\end{eqnarray*}
The left hand side of this is equal to $s({\rm dav}_{fc}(\mathbf{d}))$, and the right hand side is equal to $s(\mathbf{d})$ by (4.1). This completes the proof of Proposition 4.2. \qed\\

The main purpose of this subsection is to give a (rather lengthy) proof of the following:

\begin{prp}
We have the following equality of maps from ${\rm CA}_N^{(n)}$ to ${\rm CD}_N^{(n)}$:
\begin{eqnarray}
{\rm diff}\circ{\rm dav}_A={\rm dav}_{fc}\circ{\rm diff}.
\end{eqnarray}
\end{prp}

\noindent
{\it Proof}. For any $\mathbf{a}=(a_0,\cdots,a_{n-1})\in{\rm CA}_N^{(n)}$, the map on the left hand side of (4.3) transforms $\mathbf{a}$ to
\begin{eqnarray*}
&&({\rm diff}\circ{\rm dav}_A)(\mathbf{a})\\
&&={\rm diff}(av_{\mathbf{Z}_N}(a_0,a_1), av_{\mathbf{Z}_N}(a_1,a_2), \cdots, av_{\mathbf{Z}_N}(a_{n-1},a_0))\\
&&=(av_{\mathbf{Z}_N}(a_1,a_2)-_Nav_{\mathbf{Z}_N}(a_0,a_1), av_{\mathbf{Z}_N}(a_2,a_3)-_Nav_{\mathbf{Z}_N}(a_1,a_2), \\
&&\hspace{10mm} \cdots, av_{\mathbf{Z}_N}(a_0,a_1)-_Nav_{\mathbf{Z}_N}(a_{n-1},a_0)).
\end{eqnarray*}
On the other hand the map on the right hand side of (4.3) transforms it to
\begin{eqnarray*}
&&({\rm dav}_{D}\circ{\rm diff})(\mathbf{a})\\
&&={\rm dav}_{fc}(a_1-_Na_0, a_2-_Na_1,\cdots, a_0-_Na_{n-1})\\
&&=(av_{fc}(a_1-_Na_0,a_2-_Na_1),av_{fc}(a_2-_Na_1,a_3-_Na_2), \\
&&\hspace{10mm} \cdots, av_{fc}(a_0-_Na_{n-1},a_1-_Na_0)).
\end{eqnarray*}
Therefore we need to show that the equality
\begin{eqnarray}
av_{\mathbf{Z}_N}(b,c)-_Nav_{\mathbf{Z}_N}(a,b)=av_{fc}(b-_Na,c-_Nb)
\end{eqnarray}
holds for any consecutive three terms $a,b,c\in\mathbf{Z}_N$ which appear in $\mathbf{a}$. We reformulate (4.4) in more tractable form by reinterpreting the $\mathbf{Z}$-average map as follows:

\begin{lem}
For any integers $a,b$ we put
\begin{eqnarray*}
av_f(a,b)=\left\lfloor\frac{a+b}{2}\right\rfloor.
\end{eqnarray*}
Then for any $a,b\in \mathbf{Z}_N$, the $\mathbf{Z}$-average $av_{\mathbf{Z}_N}(a,b)$ is computed through the rule
\begin{eqnarray*}
av_{\mathbf{Z}_N}(a,b)=
\left\{
\begin{array}{ll}
av_f(a,b), & \mbox{if }a\leq b,\\
R_N(av_f(a,b+N)), & \mbox{if } a>b.\\
\end{array}
\right.
\end{eqnarray*}
\end{lem}

\noindent
{\it Proof}. When $a, b\in\mathbf{Z}_N$ and $a\leq b$, we have $-a+_Nb=-a+b$. Hence we can compute as follows:
\begin{eqnarray*}
av_{\mathbf{Z}_N}(a,b)&=&a+_N\left\lfloor\frac{-a+_Nb}{2}\right\rfloor\\
&=&a+_N\left\lfloor\frac{-a+b}{2}\right\rfloor\\
&=&R_N(a+\left\lfloor\frac{-a+b}{2}\right\rfloor)\\
&=&R_N(\left\lfloor\frac{2a-a+b}{2}\right\rfloor)\\
&&(\Leftarrow\mbox{since }x+\left\lfloor\frac{y}{2}\right\rfloor=\left\lfloor\frac{2x+y}{2}\right\rfloor\mbox{ holds for any integers }x,y)\\
&=&R_N(\left\lfloor\frac{a+b}{2}\right\rfloor)\\
&=&\left\lfloor\frac{a+b}{2}\right\rfloor\\
&=&av_f(a,b)
\end{eqnarray*}
On the other hand, when $a>b$, we compute as follows:
\begin{eqnarray*}
av_{\mathbf{Z}_N}(a,b)&=&a+_N\left\lfloor\frac{-a+_Nb}{2}\right\rfloor\\
&=&a+_N\left\lfloor\frac{-a+b+N}{2}\right\rfloor\\
&=&R_N(a+\left\lfloor\frac{-a+b+N}{2}\right\rfloor)\\
&=&R_N(\left\lfloor\frac{2a-a+b+N}{2}\right\rfloor)\\
&=&R_N(\left\lfloor\frac{a+b+N}{2}\right\rfloor)\\
&=&R_N(av_f(a,b+N)).
\end{eqnarray*}
This completes the proof of Lemma 4.1. \qed\\

\noindent
With the help of Lemma 4.1, the following lemma enables us to compute both sides of (4.4) directly.
\begin{lem}
In order to show the validity of $(4.4)$, we are only to show that the equality
\begin{eqnarray}
av_f(a_2,a_3)-_Nav_f(a_1,a_2)=av_{fc}(a_2-a_1,a_3-a_2)
\end{eqnarray}
holds for any triple $(a_1,a_2,a_3)\in\mathbf{Z}^3$ which satisfies the condition $0<a_2-a_1, a_3-a_2<N$.
\end{lem}

\noindent
{\it Proof}. Since $\mathbf{a}\in {\rm CA}_N^{(n)}$ has only one jumping number, the consecutive three terms $a,b,c$, appearing in $\mathbf{a}$, satisfy exactly one of the following inequalities:
\begin{eqnarray*}
(1)&&a<b<c,\\
(2)&&c<a<b,\\
(3)&&b<c<a
\end{eqnarray*}
We will examine the validity of Lemma 4.2 for each case. \\

\noindent
(1) In this case we have $0\leq a<b<c\leq N-1$, hence $b-_Na=b-a$ and $c-_Nb=c-b$. Therefore it follows from Lemma 4.1 that (4.4) is equivalent to the equality
\begin{eqnarray}
av_f(b,c)-_Nav_f(a,b)=av_{fc}(b-a,c-b).
\end{eqnarray}
Furthermore we have $0<b-a, c-b<N$, hence the lemma is proved in the case (1). \\

\noindent
(2) In this case, note that the equality $c-_Nb=(c+N)-b$ holds. The left hand side of (4.4) is computed by Lemma 4.1 as follows:
\begin{eqnarray*}
av_{\mathbf{Z}_N}(b,c)-_Nav_{\mathbf{Z}_N}(a,b)&=&R_N(av_f(b,c+N))-_Nav_f(a,b))\\
&=&av_f(b,c+N)-_Nav_f(a,b).
\end{eqnarray*}
The right hand side of (4.4) becomes
\begin{eqnarray*}
av_{fc}(b-_Na,c-_Nb)&=&av_{fc}(b-a,(c+N)-b).
\end{eqnarray*}
Hence the equality (4.4) is equivalent to the assertion that (4.6) holds for the triple $a,b,c+N\in\mathbf{Z}$ with $a<b<c+N$. Therefore the lemma is proved  in the case (2). \\

\noindent
(3) By Lemma 4.1, the left hand side of (4.4) is computed as
\begin{eqnarray*}
&&av_{\mathbf{Z}_N}(b,c)-_Nav_{\mathbf{Z}_N}(a,b)\\
&&=av_f(b,c)-_NR_N(av_f(a,b+N))\\
&&=R_N(av_f(b+N,c+N))-_NR_N(av_f(a,b+N))\\
&&=av_f(b+N,c+N)-_Nav_f(a,b+N).
\end{eqnarray*}
On the other hand the right hand side of (4.4) becomes
\begin{eqnarray*}
av_{fc}(b-_Na,c-_Nb)=av_{fc}((b+N)-a,(c+N)-(b+N)).
\end{eqnarray*}
Hence the equality (4.4) is equivalent to the assertion that (4.6) holds for the triple $a,b+N,c+N\in\mathbf{Z}$ with $a<b+N<c+N$. Furthermore the inequalities $0<(b+N)-a, (c+N)-(b+N)<N$ hold. This shows the validity of the lemma in the case (3), and completes the proof of the lemma at the same time. \qed\\

\noindent
By Lemma 4.2, we are reduced to show the validity of the following:\\

\noindent
{\it Claim}: The equality (4.5) holds for any triple $\mathbf{a}=(a_1,a_2,a_3)\in\mathbf{Z}^3$ which satisfies the condition $0<a_2-a_1, a_3-a_2<N$. \\

\noindent
We put \begin{eqnarray*}
\mathbf{a}'=(a'_1,a'_2)=(av_f(a_1,a_2), av_f(a_2,a_3)).
\end{eqnarray*}
Then we have
\begin{eqnarray*}
a'_1&=&\left\lfloor\frac{a_1+a_2}{2}\right\rfloor,\\
a'_2&=&\left\lfloor\frac{a_2+a_3}{2}\right\rfloor.
\end{eqnarray*}
If we denote the difference of $\mathbf{a}$ by $\mathbf{d}=(d_1,d_2)$, then we have
\begin{eqnarray*}
(d_1,d_2)&=&(a_2-_Na_1, a_3-_Na_2)\\
&=&(a_2-a_1,a_3-a_2)\hspace{3mm}(\Leftarrow \mbox{since } 0<a_2-a_1, a_3-a_2\leq N-1.)
\end{eqnarray*}
Furthermore let $d'=a'_2-_Na'_1$ denote the difference of $\mathbf{a}'$. Note that $a_1\leq av_f(a_1,a_2)\leq a_2$ and $a_2\leq av_f(a_2,a_3)\leq a_3$ holds by the definition of $av_f$, we have $0\leq a_1\leq a'_1\leq a'_2\leq a_3\leq N-1$. Therefore
\begin{eqnarray*}
d'=a'_2-a'_1.
\end{eqnarray*}
It follows that the left hand side of (4.5) becomes
\begin{eqnarray*}
av_f(a_2,a_3)-_Nav_f(a_1,a_2)=a'_2-_Na'_1=d',
\end{eqnarray*}
and the right hand side of (4.5) equals
\begin{eqnarray*}
av_{fc}(a_2-a_1,a_3-a_2)=av_{fc}(d_1,d_2).
\end{eqnarray*}
Therefore in order to show the validity of (4.5), it suffices to show that
\begin{eqnarray}
av_{fc}(\mathbf{d})=d'.
\end{eqnarray}
We divide the argument according to the parities of $a_1, a_2, a_3$.\\

\noindent
(1) The case when $(a_1,a_2,a_3)\equiv (0,0,0), (1,1,1)\pmod{2}$: In this case we see that
\begin{eqnarray*}
&&\mathbf{a}'=(a'_1,a'_2)=\left(\left\lfloor\frac{a_1+a_2}{2}\right\rfloor, \left\lfloor\frac{a_2+a_3}{2}\right\rfloor\right)=\left(\frac{a_1+a_2}{2}, \frac{a_2+a_3}{2}\right),\\
&&\mathbf{d}=(d_1,d_2)=(a_2-a_1, a_3-a_2),\\
&&d'=a'_2-a'_1=\frac{a_2+a_3}{2}-\frac{a_1+a_2}{2}=\frac{a_3-a_1}{2}.
\end{eqnarray*}
Since the first entry of $\mathbf{d}$ is even, we have
\begin{eqnarray*}
av_{fc}(\mathbf{d})&=&\left\lfloor\frac{(a_2-a_1)+(a_3-a_2)}{2}\right\rfloor\\
&=&\left\lfloor\frac{a_3-a_1}{2}\right\rfloor\\
&=&\frac{a_3-a_1}{2}\hspace{5mm}(\Leftarrow \mbox{since }a_3-a_1\equiv 0\pmod 2)\\
&=&d'.
\end{eqnarray*}
Thus the equality (4.7) is proved in case (1).\\

\noindent
(2) The case when $(a_1,a_2,a_3)\equiv (0,0,1), (1,1,0)\pmod{2}$: We have
\begin{eqnarray*}
&&\mathbf{a}'=\left(\left\lfloor\frac{a_1+a_2}{2}\right\rfloor, \left\lfloor\frac{a_2+a_3}{2}\right\rfloor\right)=\left(\frac{a_1+a_2}{2}, \frac{a_2+a_3-1}{2}\right),\\
&&\mathbf{d}=(a_2-a_1, a_3-a_2),\\
&&d'=\frac{a_2+a_3-1}{2}-\frac{a_1+a_2}{2}=\frac{a_3-a_1-1}{2}.
\end{eqnarray*}
Since the first entry of $\mathbf{d}$ is even, we have
\begin{eqnarray*}
av_{fc}(\mathbf{d})&=&\left\lfloor\frac{(a_2-a_1)+(a_3-a_2)}{2}\right\rfloor\\
&=&\left\lfloor\frac{a_3-a_1}{2}\right\rfloor\\
&=&\frac{a_3-a_1-1}{2}\hspace{5mm}(\Leftarrow \mbox{since }a_3-a_1\equiv 1\pmod 2)\\
&=&d'.
\end{eqnarray*}
Thus the equality (4.7) is proved in case (2).\\

\noindent
(3) The case when $(a_1,a_2,a_3)\equiv (0,1,0), (1,0,1)\pmod{2}$: We have
\begin{eqnarray*}
&&\mathbf{a}'=\left(\left\lfloor\frac{a_1+a_2}{2}\right\rfloor, \left\lfloor\frac{a_2+a_3}{2}\right\rfloor\right)=\left(\frac{a_1+a_2-1}{2}, \frac{a_2+a_3-1}{2}\right),\\
&&\mathbf{d}=(a_2-a_1, a_3-a_2),\\
&&d'=\frac{a_2+a_3-1}{2}-\frac{a_1+a_2-1}{2}=\frac{a_3-a_1}{2}.
\end{eqnarray*}
Since the first entry of $\mathbf{d}$ is odd, we have
\begin{eqnarray*}
av_{fc}(\mathbf{d})&=&\left\lceil\frac{(a_2-a_1)+(a_3-a_2)}{2}\right\rceil\\
&=&\left\lceil\frac{a_3-a_1}{2}\right\rceil\\
&=&\frac{a_3-a_1}{2}\hspace{5mm}(\Leftarrow \mbox{since }a_3-a_1\equiv 0\pmod 2)\\
&=&d'.
\end{eqnarray*}
Thus the equality (4.7) is proved in case (3).\\

\noindent
(4) The case when $(a_1,a_2,a_3)\equiv (0,1,1), (1,0,0)\pmod{2}$: We have
\begin{eqnarray*}
&&\mathbf{a}'=\left(\left\lfloor\frac{a_1+a_2}{2}\right\rfloor, \left\lfloor\frac{a_2+a_3}{2}\right\rfloor\right)=\left(\frac{a_1+a_2-1}{2}, \frac{a_2+a_3}{2}\right)\\
&&\mathbf{d}=(a_2-a_1, a_3-a_2)\\
&&d'=\frac{a_2+a_3}{2}-\frac{a_1+a_2-1}{2}=\frac{a_3-a_1+1}{2}
\end{eqnarray*}
Since the first entry of $\mathbf{d}$ is odd, we have
\begin{eqnarray*}
av_{fc}(\mathbf{d})&=&\left\lceil\frac{(a_2-a_1)+(a_3-a_2)}{2}\right\rceil\\
&=&\left\lfloor\frac{a_3-a_1}{2}\right\rfloor\\
&=&\frac{a_3-a_1+1}{2}\hspace{5mm}(\Leftarrow \mbox{since }a_3-a_1\equiv 1\pmod 2)\\
&=&d'.
\end{eqnarray*}
Thus the equality (4.7) is proved in case (4), therefore we complete the proof of (4.7) and Proposition 4.2 at the same time. \qed\\

\noindent
Thanks to Proposition 4.2, we can compute $\mathbf{d}^{(k+1)}$ directly from $\mathbf{d}^{(k)}$. The reader will be convinced that it is actually the case in Example 3.1 and 3.2. \\

\subsection{$D$-average transformation and the width}
　For any $\mathbf{d}\in\mathbf{Z}^n$, we define its {\it width} $w(\mathbf{d})$ by 
\begin{eqnarray*}
w(\mathbf{d})=\max(\mathbf{d})-\min(\mathbf{d}).
\end{eqnarray*}
From now on, our main objective is to understand the behavior of the width under the iteration of $D$-average map. Since the $fc$-average satisfies the inequality
\begin{eqnarray}
\min (a,b)\leq av_{fc}(a,b)\leq\max(a,b)
\end{eqnarray}
for any $a,b\in \mathbf{Z}$, we see that
\begin{eqnarray*}
\min(\mathbf{d})\leq\min(\mathbf{d}^{(1)})\leq\max(\mathbf{d}^{(1)})\leq\max(\mathbf{d})
\end{eqnarray*}
holds for any $\mathbf{d}\in\mathbf{Z}^n$. This provides us with the following simple but important observation:
\begin{prp}
For any $\mathbf{d}\in\mathbf{Z}^n$, we hace
\begin{eqnarray*}
w(\mathbf{d}^{(1)})\leq w(\mathbf{d}).
\end{eqnarray*}
\end{prp}

\subsection{Sequences with width 0 or 1}
In this subsection we show that the integer sequences with width 0 or 1 are fixed (or nearly fixed) by the $D$-average transformation. We need to introduce the following map in order to formulate our result:

\begin{df}
For any $\mathbf{d}=(d_0,d_1,\cdots, d_{n-1})\in\mathbf{Z}^n$, let ${\rm cyc}^+$ denote the map defined by
\begin{eqnarray*}
{\rm cyc}^+(d_0,d_1,\cdots, d_{n-1})=(d_{n-1},d_0,\cdots,d_{n-2}),
\end{eqnarray*}
and call it the {\rm right rotation}. Similarly we define the map ${\rm cyc}^-$ by
\begin{eqnarray*}
{\rm cyc}^-(d_0,d_1,\cdots, d_{n-1})=(d_1,d_2,\cdots,d_{n-1},d_0),
\end{eqnarray*}
and call it the {\rm left rotation}.
\end{df}

\noindent
The following proposition classifies the fixed points of width 0 or 1 completely:

\begin{prp}
$(1)$ Suppose that $\mathbf{d}\in\mathbf{Z}^n$ is of width $0$. Then we have
\begin{eqnarray*}
\mathbf{d}^{(1)}=\mathbf{d}.
\end{eqnarray*}
Furthermore $n$ must divide the sum $s(\mathbf{d})$ in this case.\\
$(2)$ Suppose that $\mathbf{d}\in\mathbf{Z}^n$ is of width $1$. Then we have the following two alternatives:\\
$(2.1)$ If $\min(\mathbf{d})$ is even, then
\begin{eqnarray*}
\mathbf{d}^{(1)}=\mathbf{d}.
\end{eqnarray*}
$(2.2)$ If $\min(\mathbf{d})$ is odd, then
\begin{eqnarray*}
\mathbf{d}^{(1)}={\rm cyc}^-(\mathbf{d}).
\end{eqnarray*}
\end{prp}

\noindent
{\it Proof}. (1) Since $\max(\mathbf{d})=\min(\mathbf{d})$ in this case, every entry of $\mathbf{d}$ has one and the same value $m$ (say). Therefore all the $fc$-averages of its consectuve entries are equal to $m$, and hence we have $\mathbf{d}^{(1)}=\mathbf{d}$. Furthermore the sum becomes $s(\mathbf{d})=\sum_{i=0}^{n-1}d_i=nm$, it follows that $n$ divides $s(\mathbf{d})$.\\

\noindent
(2.1) Let $\min(\mathbf{d})=m$. Since $m$ is even by the assumption, we have 
\begin{eqnarray*}
av_{fc}(m,m)&=&m,\\
av_{fc}(m,m+1)&=&m,\\
av_{fc}(m+1,m)&=&m+1,\\
av_{fc}(m+1,m+1)&=&m+1.
\end{eqnarray*}
Furthermore, since $d_i, d_{i+1}\in\{m,m+1\}$ for any $i\in [0, n-1]$, these equalities imply that
\begin{eqnarray*}
av_{fc}(d_i, d_{i+1})=d_i.
\end{eqnarray*}
Hence $d^{(1)}_i=d_i$ holds for any $i\in [0,n-1]$ and we have $\mathbf{d}^{(1)}=\mathbf{d}$.\\

\noindent
(2.2) Let $\min(\mathbf{d})=m$. Since $m$ is odd by the assumption, it follows that
\begin{eqnarray*}
av_{fc}(m,m)&=&m,\\
av_{fc}(m,m+1)&=&m+1,\\
av_{fc}(m+1,m)&=&m,\\
av_{fc}(m+1,m+1)&=&m+1.
\end{eqnarray*}
Since $d_i, d_{i+1}\in\{m,m+1\}, i\in [0,n-1]$ by the assumption, these equalities imply that \begin{eqnarray*}
av_{fc}(d_i, d_{i+1})=d_{i+1}.
\end{eqnarray*}
Therefore we have $\mathbf{d}^{(1)}={\rm cyc}^-(\mathbf{d})$. This completes the proof. \qed\\

\section{Condition $C(\mathbf{d})$}
Our main theorem will be established by showing that, if we apply the $D$-average transformaition repeatedly to any given $\mathbf{d}\in\mathbf{Z}^n$, then the width of $\mathbf{d}^{(k)}$ becomes smaller than or equal to one for an appropriate $k$. \\

\subsection{Formulation of $C(\mathbf{d})$}
　We introduce some terminologies and formulate an inductive step, called $C(\mathbf{d})$ from now on, in order to accomplish the proof of the main theorem. \\

For any integer $n\geq 3$, let $C_n=(V(C_n),E(C_n))$ denote the graph with vertex set $V(C_n)=\mathbf{Z}_n$. The set of edges $E(C_n)$ consist of the pairs $(i,j)\in\mathbf{Z}_n^2$ with $|i-j|=1$. Thus $C_n$ is isomorphic to the graph which consists of the vertices and the edges of the regular $n$-gon. For any pairs $(i,j)$ of the vertices of $C_n$, let $[i,j]$ denote the connected subset of $V(C_n)$ defined by the rule
\begin{eqnarray*}
[i,j]=
\left\{
\begin{array}{ll}
\{i,i+1,\cdots,j-1,j\}, & \mbox{if }0\leq i\leq j\leq n-1,\\
\{i,i+1,\cdots,n-1,0,1,\cdots,j\}, & \mbox{if }0\leq j<i\leq n-1.\\
\end{array}
\right.
\end{eqnarray*}
We fix some notation:

\begin{df}
For any $\mathbf{d}=(d_0,\cdots,d_{n-1})\in \mathbf{Z}^n$, we put
\begin{eqnarray*}
\mathbf{I}_{min}(\mathbf{d})&=&\{k\in C_n;d_k=\min(\mathbf{d})\},\\
\mathbf{I}_{max}(\mathbf{d})&=&\{k\in C_n;d_k=\max(\mathbf{d})\}.
\end{eqnarray*}
Furthermore regarding the set $\mathbf{I}_{min}(\mathbf{d})$ as the subgraph of $C_n$ we represent it as a disjoint union of its connected components:
\begin{eqnarray*}
\mathbf{I}_{min}(\mathbf{d})=\bigsqcup_{p=1}^{n_{min}}S_p.
\end{eqnarray*}
We call each $S_p, 1\leq p\leq n_{min},$ a {\rm block of minimum values}. Similarly we represent the set $\mathbf{I}_{max}(\mathbf{d})$ as a disjoint union of its connected components:
\begin{eqnarray*}
\mathbf{I}_{max}(\mathbf{d})=\bigsqcup_{q=1}^{n_{max}}L_q,
\end{eqnarray*}
and call each component $L_q, 1\leq q\leq n_{max},$ a {\rm block of maximum values}. Moreover we put
\begin{eqnarray*}
s_{min}(\mathbf{d})&=&\sum_{p=1}^{n_{min}}|S_p|,\\
s_{max}(\mathbf{d})&=&\sum_{q=1}^{n_{max}}|L_q|.
\end{eqnarray*}
\end{df}

We will prove our main theorem by a double induction on the width and the sum of length $s_{min}(\mathbf{d})$ of the blocks of minimum values. More precisely we focus on the following condition $C(\mathbf{d})$ on $\mathbf{d}\in\mathbf{Z}^n$:\\

\noindent
$C(\mathbf{d})$: There exists a positive integer $k$ such that either $w(\mathbf{d}^{(k)})<w(\mathbf{d})$ or $s_{min}(\mathbf{d}^{(k)})<s_{min}(\mathbf{d})$ holds.\\

\subsection{Reduction of proof of $C(\mathbf{d})$ for odd cases to even cases}
When we try to prove the condition $C(\mathbf{d})$, however, the behavior of $\mathbf{d}^{(k)}, k\geq 1,$ presents an entirely different picture depending on the parity of $\min(\mathbf{d})$, as is imagined from the conclusion of Proposition 4.4, (2). In view of this, we introduce some maps in order to make our argument valid notwithstanding the parity of $\min(\mathbf{d})$: 

\begin{df}
For any $\mathbf{d}=(d_0,d_1,\cdots,d_{n-2},d_{n-1})\in\mathbf{Z}^n$, we define a map $\iota:\mathbf{Z}^n\rightarrow\mathbf{Z}^n$, which reverse the order of its entries, as follows:
\begin{eqnarray*}
\iota(d_0,d_1,\cdots,d_{n-2},d_{n-1})=(d_{n-1},d_{n-2},\cdots,d_1,d_0).
\end{eqnarray*}
The rule of numbering of its value is given by
\begin{eqnarray*}
(\iota(\mathbf{d}))_i=d_{(n-1)-i},\hspace{2mm}0\leq i\leq n-1.
\end{eqnarray*}
Furthermore, by exchanging the parities of the $fc$-average, we define the $cf$-average $av_{cf}:\mathbf{Z}^2\rightarrow\mathbf{Z}$ by
\begin{eqnarray*}
av_{cf}(a,b)=
\left\{
\begin{array}{ll}
\lceil\frac{a+b}{2}\rceil, & \mbox{if }a\equiv 0\pmod{2},\\
\lfloor\frac{a+b}{2}\rfloor, & \mbox{if }a\equiv 1\pmod{2}.\\
\end{array}
\right.
\end{eqnarray*}
Accordingly we define another $D$-average transformation ${\rm dav}_{cf}:\mathbf{Z}^n\rightarrow\mathbf{Z}^n$ by
\begin{eqnarray*}
{\rm dav}_{cf}(\mathbf{d})=(av_{cf}(d_0,d_1), av_{cf}(d_1,d_2), \cdots, av_{cf}(d_{n-1},d_0)).
\end{eqnarray*}
\end{df}

\noindent
The following proposition clarifies the relation between this average ${\rm dav}_{cf}$ and the former average ${\rm dav}_{fc}$ introduced in Definition 4.2: 

\begin{prp}
As maps from $\mathbf{Z}^n$ to $\mathbf{Z}^n$, the equality
\begin{eqnarray}
\iota\circ{\rm cyc}^+\circ{\rm dav}_{fc}={\rm dav}_{cf}\circ\iota
\end{eqnarray}
holds. In another word, we have the following commutative diagram:

\includegraphics{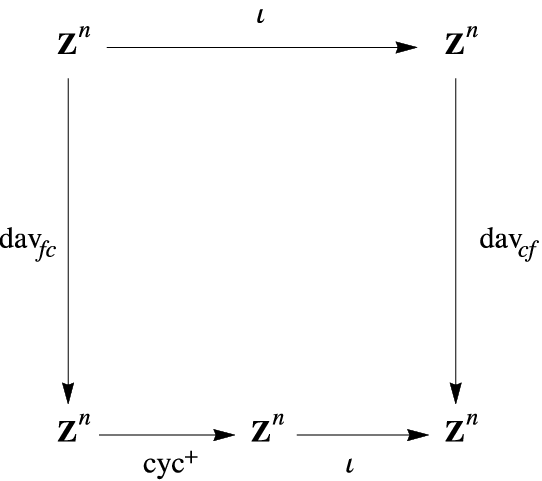}

\end{prp}

\noindent
{\it Proof}. Focusing on the parities of the consecutive pair $(d_i,d_{i+1}), i\in [0,n-1]$, of entries of $\mathbf{d}\in\mathbf{Z}^n$, we trace how the maps on both sides of (5.1) transform $\mathbf{d}$. \\

\noindent
{\it Case} (1): $i<n-1$.\\
(1.1) The case when $(d_i,d_{i+1})\equiv (0,0),(1,1)\pmod 2$. The map on the left hand side of (5.1) transforms $\mathbf{d}$ as follows:
\begin{eqnarray*}
\begin{array}{c}
(d_0,\cdots,\overbrace{d_i}^{i},d_{i+1},\cdots,d_{n-1})\\
\downarrow\hspace{2mm}{\rm dav}_{fc}\\
(*,\cdots,\overbrace{\frac{d_i+d_{i+1}}{2}}^{i},*,\cdots,*)\\
\downarrow\hspace{2mm}{\rm cyc}^+\\
(*,\cdots,\overbrace{\frac{d_i+d_{i+1}}{2}}^{i+1},*,\cdots,*)\\
\downarrow\hspace{2mm}\iota\\
(*,\cdots,\overbrace{\frac{d_i+d_{i+1}}{2}}^{(n-1)-(i+1)},*,\cdots,*).\\
\end{array}
\end{eqnarray*}
On the other hand, the map on the right hand side of (5.1) transforms it as
\begin{eqnarray*}
\begin{array}{c}
(d_0,\cdots,\overbrace{d_i}^{i},d_{i+1},\cdots,d_{n-1})\\
\downarrow\hspace{2mm}\iota\\
(*,\cdots,\overbrace{d_{i+1}}^{(n-1)-(i+1)},\overbrace{d_{i}}^{(n-1)-i},*,\cdots,*)\\
\downarrow\hspace{2mm}{\rm dav}_{cf}\\
(*,\cdots,\overbrace{\frac{d_i+d_{i+1}}{2}}^{(n-1)-(i+1)},*,\cdots,*).\\
\end{array}
\end{eqnarray*}
Therefore the value of both sides  at $\mathbf{d}$ coincide at the $i$-th place.\\

\noindent
(1.2) The case when $(d_i,d_{i+1})\equiv (0,1)\pmod 2$. The map on the left hand side of (5.1) transforms $\mathbf{d}$ as follows:
\begin{eqnarray*}
\begin{array}{c}
(d_0,\cdots,\overbrace{d_i}^{i},d_{i+1},\cdots,d_{n-1})\\
\downarrow\hspace{2mm}{\rm dav}_{fc}\\
(*,\cdots,\overbrace{\frac{d_i+d_{i+1}-1}{2}}^{i},*,\cdots,*)\\
\downarrow\hspace{2mm}{\rm cyc}^+\\
(*,\cdots,\overbrace{\frac{d_i+d_{i+1}-1}{2}}^{i+1},*,\cdots,*)\\
\downarrow\hspace{2mm}\iota\\
(*,\cdots,\overbrace{\frac{d_i+d_{i+1}-1}{2}}^{(n-1)-(i+1)},*,\cdots,*).\\
\end{array}
\end{eqnarray*}
On the other hand, the map on the right hand side of (5.1) transforms it as
\begin{eqnarray*}
\begin{array}{c}
(d_0,\cdots,\overbrace{d_i}^{i},d_{i+1},\cdots,d_{n-1})\\
\downarrow\hspace{2mm}\iota\\
(*,\cdots,\overbrace{d_{i+1}}^{(n-1)-(i+1)},\overbrace{d_{i}}^{(n-1)-i},*,\cdots,*)\\
\downarrow\hspace{2mm}{\rm dav}_{cf}\\
(*,\cdots,\overbrace{\frac{d_i+d_{i+1}-1}{2}}^{(n-1)-(i+1)},*,\cdots,*).\\
\end{array}
\end{eqnarray*}
Therefore the value ot both sides  at $\mathbf{d}$ coincide at the $i$-th place.\\

\noindent
(1.3) The case when $(d_i,d_{i+1})\equiv (1,0)\pmod 2$. The map on the left hand side of (5.1) transforms $\mathbf{d}$ as follows:
\begin{eqnarray*}
\begin{array}{c}
(d_0,\cdots,\overbrace{d_i}^{i},d_{i+1},\cdots,d_{n-1})\\
\downarrow\hspace{2mm}{\rm dav}_{fc}\\
(*,\cdots,\overbrace{\frac{d_i+d_{i+1}+1}{2}}^{i},*,\cdots,*)\\
\downarrow\hspace{2mm}{\rm cyc}^+\\
(*,\cdots,\overbrace{\frac{d_i+d_{i+1}+1}{2}}^{i+1},*,\cdots,*)\\
\downarrow\hspace{2mm}\iota\\
(*,\cdots,\overbrace{\frac{d_i+d_{i+1}+1}{2}}^{(n-1)-(i+1)},*,\cdots,*).\\
\end{array}
\end{eqnarray*}
On the other hand, the map on the right hand side of (5.1) transforms it as
\begin{eqnarray*}
\begin{array}{c}
(d_0,\cdots,\overbrace{d_i}^{i},d_{i+1},\cdots,d_{n-1})\\
\downarrow\hspace{2mm}\iota\\
(*,\cdots,\overbrace{d_{i+1}}^{(n-1)-(i+1)},\overbrace{d_{i}}^{(n-1)-i},*,\cdots,*)\\
\downarrow\hspace{2mm}{\rm dav}_{cf}\\
(*,\cdots,\overbrace{\frac{d_i+d_{i+1}+1}{2}}^{(n-1)-(i+1)},*,\cdots,*)\\
\end{array}
\end{eqnarray*}
Therefore the value ot both sides  at $\mathbf{d}$ coincide at the $i$-th place.\\

\noindent
{\it Case} (2): $i=n-1$.\\
(2.1) The case when $(d_i,d_{i+1})\equiv (0,0),(1,1)\pmod 2$. The map on the left hand side of (5.1) transforms $\mathbf{d}$ as follows:
\begin{eqnarray*}
\begin{array}{c}
(d_0,\cdots,d_{n-1})\\
\downarrow\hspace{2mm}{\rm dav}_{fc}\\
(*,\cdots,\frac{d_{n-1}+d_0}{2})\\
\downarrow\hspace{2mm}{\rm cyc}^+\\
(\frac{d_{n-1}+d_0}{2},\cdots,*)\\
\downarrow\hspace{2mm}\iota\\
(*,\cdots,\frac{d_{n-1}+d_0}{2}).\\
\end{array}
\end{eqnarray*}
On the other hand, the map on the right hand side of (5.1) transforms it as
\begin{eqnarray*}
\begin{array}{c}
(d_0,\cdots,d_{n-1})\\
\downarrow\hspace{2mm}\iota\\
(d_{n-1},\cdots,d_0)\\
\downarrow\hspace{2mm}{\rm dav}_{cf}\\
(*,\cdots,\frac{d_0+d_{n-1}}{2}).\\
\end{array}
\end{eqnarray*}
Therefore the value ot both sides  at $\mathbf{d}$ coincide at the $(n-1)$-th place.\\

\noindent
(2.2) The case when $(d_{n-1},d_0)\equiv (0,1)\pmod 2$. The map on the left hand side of (5.1) transforms $\mathbf{d}$ as follows:
\begin{eqnarray*}
\begin{array}{c}
(d_0,\cdots,d_{n-1})\\
\downarrow\hspace{2mm}{\rm dav}_{fc}\\
(*,\cdots,\frac{d_{n-1}+d_0-1}{2})\\
\downarrow\hspace{2mm}{\rm cyc}^+\\
(\frac{d_{n-1}+d_0-1}{2},\cdots,*)\\
\downarrow\hspace{2mm}\iota\\
(*,\cdots,\frac{d_{n-1}+d_0-1}{2}).\\
\end{array}
\end{eqnarray*}
On the other hand, the map on the right hand side of (5.1) transforms it as
\begin{eqnarray*}
\begin{array}{c}
(d_0,\cdots,d_{n-1})\\
\downarrow\hspace{2mm}\iota\\
(d_{n-1},\cdots,d_0)\\
\downarrow\hspace{2mm}{\rm dav}_{cf}\\
(*,\cdots,\frac{d_0+d_{n-1}-1}{2}).\\
\end{array}
\end{eqnarray*}
Therefore the value ot both sides  at $\mathbf{d}$ coincide at the $(n-1)$-th place.\\

\noindent
(2.3) The case when $(d_{n-1},d_0)\equiv (1,0)\pmod 2$. The map on the left hand side of (5.1) transforms $\mathbf{d}$ as follows:
\begin{eqnarray*}
\begin{array}{c}
(d_0,\cdots,d_{n-1})\\
\downarrow\hspace{2mm}{\rm dav}_{fc}\\
(*,\cdots,\frac{d_{n-1}+d_0+1}{2})\\
\downarrow\hspace{2mm}{\rm cyc}^+\\
(\frac{d_{n-1}+d_0+1}{2},\cdots,*)\\
\downarrow\hspace{2mm}\iota\\
(*,\cdots,\frac{d_{n-1}+d_0+1}{2}).\\
\end{array}
\end{eqnarray*}
On the other hand, the map on the right hand side of (5.1) transforms it as
\begin{eqnarray*}
\begin{array}{c}
(d_0,\cdots,d_{n-1})\\
\downarrow\hspace{2mm}\iota\\
(d_{n-1},\cdots,d_0)\\
\downarrow\hspace{2mm}{\rm dav}_{cf}\\
(*,\cdots,\frac{d_0+d_{n-1}+1}{2}).\\
\end{array}
\end{eqnarray*}
Therefore the value ot both sides  at $\mathbf{d}$ coincide at the $(n-1)$-th place. This completes the proof of Proposition 5.1.\qed\\

Furthermore the following proposition relates the two maps ${\rm dav}_{fc}$ and ${\rm dav}_{cf}$ directly:

\begin{prp}
For any pair of integers $p, q$, we have
\begin{eqnarray}
av_{fc}(p+1,q+1)=av_{cf}(p,q)+1.
\end{eqnarray}
Furthermore, if we define a map ${\rm add}^+:\mathbf{Z}^n\rightarrow\mathbf{Z}^n$ by
\begin{eqnarray*}
{\rm add}^+(d_0,\cdots,d_{n-1})=(d_0+1,\cdots,d_{n-1}+1),
\end{eqnarray*}
then the equality
\begin{eqnarray}
{\rm add}^+\circ{\rm dav}_{cf}={\rm dav}_{fc}\circ{\rm add}^+
\end{eqnarray}
holds. Namely we have the following commutative diagram:

\includegraphics{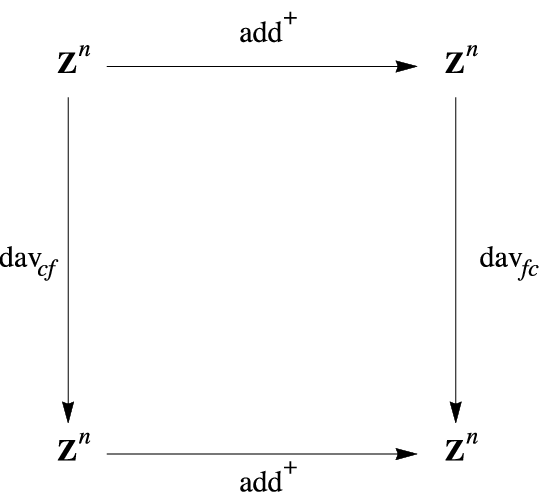}

\end{prp}

\noindent
{\it Proof}. When $p+1$ is even, the left hand side of (5.2) becomes
\begin{eqnarray*}
av_{fc}(p+1,q+1)&=&\left\lfloor\frac{(p+1)+(q+1)}{2}\right\rfloor\\
&=&\left\lfloor\frac{p+q}{2}\right\rfloor+1.
\end{eqnarray*}
On the other hand, the right hand side of (5.2) becomes
\begin{eqnarray*}
av_{cf}(p,q)+1=\left\lfloor\frac{p+q}{2}\right\rfloor+1,
\end{eqnarray*}
since $p$ is odd. Hence the equality (5.2) holds in this case.\\

 \noindent
Changing the floor functions, which appear in the above proof for the case $p+1$ is even, to the ceiling functions, we obtain automatically a proof for the case when $p+1$ is odd. Thus the equality (5.2) holds in both cases. For the proof of the equality (5.3), we apply the right hand side to an arbitrary $\mathbf{d}=(d_0,\cdots,d_{n-1})\in\mathbf{Z}^n$. Then we have
\begin{eqnarray*}
&&{\rm dav}_{fc}({\rm add}^+(d_0,\cdots,d_{n-1}))\\
&=&{\rm dav}_{fc}(d_0+1,\cdots,d_{n-1}+1)\\
&=&(av_{fc}(d_0+1,d_1+1),\cdots,av_{fc}(d_{n-1}+1,d_0+1))\\
&=&(av_{cf}(d_0,d_1)+1,\cdots,av_{cf}(d_{n-1},d_0)+1)\\
&&\hspace{30mm}(\Leftarrow\mbox{ by }(5.2))\\
&=&{\rm add}^+(av_{cf}(d_0,d_1),\cdots,av_{cf}(d_{n-1},d_0))\\
&=&{\rm add}^+({\rm dav}_{cf}(d_0,\cdots, d_{n-1}).
\end{eqnarray*}
This coincides with the value of the left hand side at $\mathbf{d}$. This completes the proof of Proposition 5.2. \qed\\

Note that the three maps $\iota, {\rm cyc}^+,$ and ${\rm add}^+$ are bijections from $\mathbf{Z}^n$ onto $\mathbf{Z}^n$, since we have
\begin{eqnarray*}
\iota\circ\iota&=&{\rm id}_{\mathbf{Z}^n},\\
{\rm cyc}^+\circ{\rm cyc}^-&=&{\rm cyc}^-\circ{\rm cyc}^+={\rm id}_{\mathbf{Z}^n},\\
{\rm add}^+\circ{\rm add}^-&=&{\rm add}^-\circ{\rm add}^+={\rm id}_{\mathbf{Z}^n}.
\end{eqnarray*}
Here the map ${\rm add}^-:\mathbf{Z}^n\rightarrow\mathbf{Z}^n$ is defined by
\begin{eqnarray*}
{\rm add}^-(d_0,\cdots,d_{n-1})=(d_0-1,\cdots,d_{n-1}-1).
\end{eqnarray*}
It follows from the definitions of these three maps that the width of $\mathbf{d}\in\mathbf{Z}^n$ is invariant under these maps:
\begin{eqnarray*}
w(\iota(\mathbf{d}))&=&w(\mathbf{d}),\\
w({\rm cyc}^+(\mathbf{d}))&=&w(\mathbf{d}),\\
w({\rm add}^+(\mathbf{d}))&=&w(\mathbf{d}).
\end{eqnarray*}
It also follows that the length of the minimum values is invariant:
\begin{eqnarray*}
\ell(\iota(\mathbf{d}))&=&\ell(\mathbf{d}),\\
\ell({\rm cyc}^+(\mathbf{d}))&=&\ell(\mathbf{d}),\\
\ell({\rm add}^+(\mathbf{d}))&=&\ell(\mathbf{d}).
\end{eqnarray*}
The minimum value itself varies under the map ${\rm add}^+$ only:
\begin{eqnarray*}
\min(\iota(\mathbf{d}))&=&\min(\mathbf{d}),\\
\min({\rm cyc}^+(\mathbf{d}))&=&\min(\mathbf{d}),\\
\min({\rm add}^+(\mathbf{d}))&=&\min(\mathbf{d})+1.
\end{eqnarray*}

\noindent
The above consideration shows the validity of the former half of the following proposition:

\begin{prp}
Let $f^+:\mathbf{Z}^n\rightarrow\mathbf{Z}^n$ and $g^+:\mathbf{Z}^n\rightarrow\mathbf{Z}^n$ denote the map defined by the rule
\begin{eqnarray*}
f^+&=&{\rm add}^+\circ\iota,\\
g^+&=&{\rm add}^+\circ\iota\circ{\rm cyc}^+.
\end{eqnarray*}
Then for any $\mathbf{d}\in\mathbf{Z}^n$, we have
\begin{eqnarray}
w(f^+(\mathbf{d}))&=&w(g^+(\mathbf{d}))=w(\mathbf{d}),\\
\ell(f^+(\mathbf{d}))&=&\ell(g^+(\mathbf{d}))=\ell(\mathbf{d}),\\
\min(f^+(\mathbf{d}))&=&\min(g^+(\mathbf{d}))=\min(\mathbf{d})+1.
\end{eqnarray}
Furthremore the equality
\begin{eqnarray}
{\rm dav}_{fc}\circ f^+=g^+\circ{\rm dav}_{fc}
\end{eqnarray}
holds. Namely we have the following commutative diagram:

\includegraphics{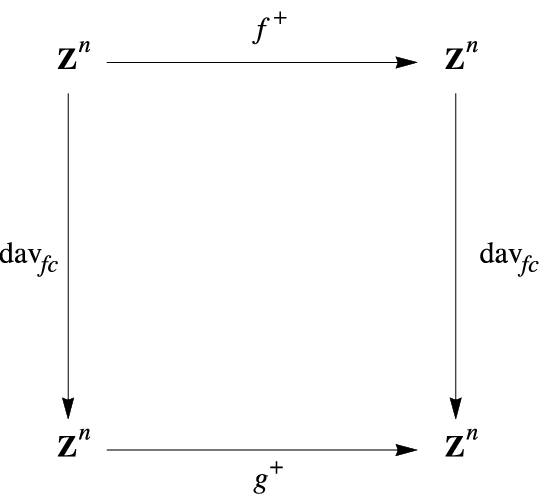}

\end{prp}

\noindent
{\it Proof}. We are only to show the validity of (5.7). The left hand side can be computed as follows:
\begin{eqnarray*}
{\rm dav}_{fc}\circ f^+&=&{\rm dav}_{fc}\circ ({\rm add}^+\circ\iota)\\
&=&({\rm dav}_{fc}\circ {\rm add}^+)\circ\iota\\
&=&({\rm add}^+\circ {\rm dav}_{cf})\circ\iota\\
&&\hspace{30mm}(\Leftarrow\mbox{ by }(5.3))\\
&=&{\rm add}^+\circ ({\rm dav}_{cf}\circ\iota)\\
&=&{\rm add}^+\circ (\iota\circ{\rm cyc}^+\circ{\rm dav}_{fc})\\
&&\hspace{30mm}(\Leftarrow\mbox{ by }(5.1))\\
&=&({\rm add}^+\circ \iota\circ{\rm cyc}^+)\circ{\rm dav}_{fc}\\
&=&g^+\circ{\rm dav}_{fc}.
\end{eqnarray*}
The rightmost side coincides with the right hand side of (5.7), and hence we finish the proof. \qed\\

\noindent
Furthermore, by applying the equality (5.7) iteratively, we obtain the following:

\begin{cor}
For any positive integer $k$, we have
\begin{eqnarray}
{\rm dav}^k_{fc}\circ f^+=g^+\circ{\rm dav}^k_{fc}.
\end{eqnarray}
\end{cor}

\noindent
　By employing Proposition 5.3 and Corollary 5.1, we will see that it suffices to show the validity of the condition $C(\mathbf{d})$ only when $\min(\mathbf{d})$ is even:

\begin{prp}
If the condition $C(\mathbf{d})$ holds for any $\mathbf{d}\in\mathbf{Z}^n$ such that $\min(\mathbf{d})$ is even, then it holds for every $\mathbf{d}\in\mathbf{Z}^n$ without any assumption about the parity of $\min(\mathbf{d})$.
\end{prp}

\noindent
{\it Proof}. Suppose that $\min(\mathbf{d})=m$ is odd. If we put $f^+(\mathbf{d})=\mathbf{d}'$, then it follows from (5.6) that
\begin{eqnarray*}
\min(\mathbf{d}')=\min(\mathbf{d})+1=m+1.
\end{eqnarray*}
and hence it is even. By the assumption that $C(\mathbf{d})$ holds for even cases, we see that there exists a positive integer $k$ such that 
\begin{eqnarray}
w(\mathbf{d}'^{(k)})<w(\mathbf{d}')\mbox{ or }\ell(\mathbf{d}'^{(k)})<\ell(\mathbf{d}')\mbox{ holds.}
\end{eqnarray}
Since we have
\begin{eqnarray*}
\mathbf{d}'^{(k)}&=&{\rm dav}^k_{fc}(f^+(\mathbf{d}))\\
&=&({\rm dav}^k_{fc}\circ f^+)(\mathbf{d}))\\
&=&(g^+\circ {\rm dav}^k_{fc})(\mathbf{d}))\\
&&\hspace{30mm}(\Leftarrow\mbox{ by }(5.8))\\
&=&g^+({\rm dav}^k_{fc}(\mathbf{d}))\\
&=&g^+(\mathbf{d}^{(k)}),
\end{eqnarray*}
the condition (5.9) can be restated as
\begin{eqnarray*}
w(g^+(\mathbf{d}^{(k)}))<w(f^+(\mathbf{d}))\mbox{ or }\ell(g^+(\mathbf{d}^{(k)}))<\ell(f^+(\mathbf{d})).
\end{eqnarray*}
Furthermore it follows from (5.4) and (5.5) that the maps $f^+, g^+$ leave the width and the minimum invariant. Hence we have 
\begin{eqnarray*}
w(\mathbf{d}^{(k)})<w(\mathbf{d})\mbox{ or }\ell(\mathbf{d}^{(k)})<\ell(\mathbf{d}),
\end{eqnarray*}
which shows the validity of the condition $C(\mathbf{d})$ for the odd cases too. This completes the proof. \qed\\

\noindent
Thus we are reduced to show the validity of the condition $C(\mathbf{d})$ for the even cases. This will be accomplished in the next section. \\

\section{Proof of $C(\mathbf{d})$ when $\min(\mathbf{d})$ is even}
Throughout this section we assume that $\mathbf{d}$ satisfies the following two conditions (A) and (B):\\

\noindent
(A): $\min(\mathbf{d})$ is even.\\
(B): $w(\mathbf{d})\geq 2$.\\

\noindent
The following proposition plays a crucial role for our argument:
\begin{prp}
For any $\mathbf{d}=(d_0,\cdots,d_{n-1})\in \mathbf{Z}^n$ with nonzero width, let $S_0=[i_0, j_0]$ be an arbitrary connected component of $\mathbf{I}_{min}(\mathbf{d})$, and put ${\rm dav}_{D}(\mathbf{d})=\mathbf{d}^{(1)}=(d^{(1)}_0,\cdots,d^{(1)}_{n-1})$. Then the following four conditions hold for $\mathbf{d}^{(1)}$:
\begin{eqnarray*}
&&(1)\hspace{2mm}d^{(1)}_{i_0-1}>m,\\
&&(2)\hspace{2mm}d^{(1)}_{j_0+1}>m,\\
&&(3)\hspace{2mm}m=d^{(1)}_{i_0}=d^{(1)}_{i_0+1}=\cdots = d^{(1)}_{j_0-1},\\
&&(4)\hspace{2mm}d^{(1)}_{j_0-1}\leq d^{(1)}_{j_0},\\
&&\hspace{10mm}(\mbox{where the equality holds if and only if }d_{j_0+1}=m+1.)
\end{eqnarray*}
In particular, except when $d_{j_0+1}=m+1$, the block $S_0$ of minimum values is transformed by ${\rm dav}_{fc}$ to $S^{(1)}_0=(d^{(1)}_{i_0},\cdots, d^{(1)}_{j_0-1})$, and its length diminishes by one.
\end{prp}

\noindent
{\it Proof}. (1) Since $S_0$ is one of blocks of minimum value, we have $d_{i_0-1}>m$. If $d_{i_0-1}>m+1$, then $d^{(1)}_{i_0-1}=av_{fc}(d_{i_0-1},d_{i_0})\geq av_{fc}(m+2,m)=\lfloor \frac{(m+2)+m}{2}\rfloor =m+1$ by the definition of $av_{fc}$, and hence the condition (1) holds. On the other hand, if $d_{i_0-1}=m+1$, then noting that $m+1$ is odd, we have $d^{(1)}_{i_0-1}=av_{fc}(d_{i_0-1},d_{i_0})=av_{fc}(m+1,m)=\lceil \frac{(m+1)+m}{2}\rceil =m+1$, and hence (1) holds too. \\
(2) By the definition of $S_0$, we have $d_{j_0+1}>m$. If $d_{j_0+1}>m+1$, then noting that $m+2$ is even, we have $d^{(1)}_{j_0+1}=av_{fc}(d_{j_0+1},d_{j_0+2})\geq av_{fc}(m+2,m)=\lfloor \frac{(m+2)+m}{2}\rfloor =m+1$, which implies the validity of the condition (2). On the other hand, if $d_{j_0+1}=m+1$, then noting that $m+1$ is odd, we have $d^{(1)}_{j_0+1}=av_{fc}(d_{j_0+1},d_{j_0+2})\geq av_{fc}(m+1,m)=\lceil \frac{(m+1)+m}{2}\rceil =m+1$, which also implies (2).\\
(3) This is simply because ${\rm av}_{fc}(m,m)=m$. \\
(4) When $d_{j_0+1}>m+1$, the definition of ${\rm av}_{fc}$ implies that $d^{(1)}_{j_0}\geq m+1$, which shows that the inequality $d^{(1)}_{j_0-1}<d^{(1)}_{j_0}$ holds in the condition (4). On the other hand, when $d_{j_0+1}=m+1$, noting that $m$ is even, we have $d^{(1)}_{j_0}={\rm av}_{fc}(d_{j_0}, d_{j_0+1})={\rm av}_{fc}(m, m+1)=\lfloor \frac{m+(m+1)}{2}\rfloor =m$, which shows that the equality $d^{(1)}_{j_0-1}=d^{(1)}_{j_0}$ holds in the condition (4). This completes the proof of Proposition 6.1.\
\qed\\

Decompose $\mathbf{I}_{min}(\mathbf{d})$ into the connected components and put
\begin{eqnarray*}
\mathbf{I}_{min}(\mathbf{d})=\bigsqcup_{p=1}^{n_{min}}S_p.
\end{eqnarray*}
We express $S_p, 1\leq p\leq n_{min},$ as
\begin{eqnarray*}
S_p=[i_p,j_p].
\end{eqnarray*}
Furthermore let $g_p$ denote the length of the gap between the consecutive components $S_p$ and $S_{p+1}$, namely we put
\begin{eqnarray*}
g_p=|[j_p+1,i_{p+1}-1]|,\hspace{3mm}p\in [0,n-1].
\end{eqnarray*}
If there is no $d_j$ between $S_p$ and $S_{p+1}$ such that $d_j\geq m+2$, namely if $d_j=m+1$ for any $j\in[j_p+1,i_{p+1}-1]$, it follows from Proposition 6.1 that 
\begin{eqnarray*}
&&d^{(1)}_{i_p}=d^{(1)}_{i_p+1}=\cdots =d^{(1)}_{j_p}=m,\\
&&d^{(1)}_{j_p+1}=\cdots =d^{(1)}_{i_{p+1}-1}=m+1.
\end{eqnarray*}
In other words, even if we apply ${\rm dav}_{fc}$ once to $\mathbf{d}$, the values $d_j, \in [i_p,i_{p+1}-1], $ do not change. Therefore we need to focus on the indices $j$ with $d_j\geq m+2$. There exists, however, at least one index $j$ such that $d_{j}\geq m+2$, since we have assumed that the width of $\mathbf{d}$ is greater than or equal to two. In view of this fact, we let 
\begin{eqnarray*}
Q=\{q_1,\cdots,q_m\}
\end{eqnarray*}
the subset of indices $p$ such that the interval $[j_p+1, i_{p+1}-1]$ contains a $d_j$ which is geater than or equal to $m+2$. It follows that for each $q_k, 1\leq k\leq m,$ there exists a number $\ell_k\in [0,g_{q_k}-1]$ such that 
\begin{eqnarray*}
&&d_{j_{q_k}+1}=\dots=d_{j_{q_k}+\ell_k}=m+1,\\
&&d_{j_{q_k}+\ell_k+1}\geq m+2.
\end{eqnarray*}
Here note that if the equality $\ell_k=0$ holds, then it means simply that $d_{j_{q_k}+1}=m+2$. In other words, the number $\ell_k$ signifies how many indices are there after $j_{q_k}+1$ for which $d_j$ remains to have one and the same value $m+1$.

The following lemma plays a crucial role in our argument later:

\begin{lem}
For any $q_k\in Q$ and for the corresponding block $S_{q_k}=[i_{q_k}, j_{q_k}]$, we have
\begin{eqnarray*}
&&d^{(\ell_k+1)}_{i_{q_k}}=\cdots=d^{(\ell_k+1)}_{j_{q_k}-1}=m, \\
&&d^{(\ell_k+1)}_{j_{q_k}}>m.
\end{eqnarray*}
\end{lem}

\noindent
{\it Proof}. By Proposition * and the fact that $m+1$ is odd, we see that if we apply ${\rm dav}_{fc}$ consecutively $ \ell_k$ times, then we have 
\begin{eqnarray*}
&&d^{(\ell_k)}_{i_{q_k}}=d^{(\ell_k)}_{i_{q_k}+1}=\cdots=d^{(\ell_k)}_{j_{q_k}}=m, \\
&&d^{(\ell_k)}_{j_{q_k}+1}\geq m+2.
\end{eqnarray*}
Therefore if we apply ${\rm dav}_{fc}$ one more time, then we have
\begin{eqnarray*}
&&d^{(\ell_k+1)}_{i_{q_k}}=\cdots=d^{(\ell_k+1)}_{j_{q_k}-1}=m, \\
&&d^{(\ell_k+1)}_{j_{q_k}}>m.
\end{eqnarray*}
This completes the proof of Lemma 6.1. \qed\\

The following lemma follows directly from this lemma:

\begin{lem}
For any $q_k\in Q$ and for the corresponding block $S_{q_k}=[i_{q_k}, j_{q_k}]$, we have the following two possibilities:\\
$(1)$ When $|S_{q_k}|>1$, namely when $i_{q_k}<j_{q_k}$, if we apply ${\rm dav}_{fc}$ consecutively $\ell_k+1$ times, then the length of the block $S^{(\ell_k+1)}_{q_k}$ of minimum value decreases by one, compare with the original block $S_{q_k}$. \\
$(2)$ When $|S_{q_k}|=1$, namely when $i_{q_k}=j_{q_k}$, if we apply ${\rm dav}_{fc}$ consecutively $\ell_k+1$ times, then the block $S^{(\ell_k+1)}_{q_k}$ of minimum value becomes the empty set.
\end{lem}

Now we begin the proof of the following proposition which is the main focus of this subsection:

\begin{prp}
Assume that the two conditions $({\rm A})$ and $({\rm B})$ are satisfied for $\mathbf{d}$. Then the condition $C(\mathbf{d})$ holds.
\end{prp}

\noindent
{\it Proof}. We divide our argument into the following three cases:\\

\noindent
(1) For any $k\in [1,m]$, it holds that $|S_{q_k}|=1$ and $\ell_k=0$.\\
(2) For any $k\in [1,m]$, it holds that $|S_{q_k}|=1$, but for some $k\in [1,m]$ the inequality $\ell_k\geq 1$ holds.\\
(3) For some $k\in [1,m]$, it holds that $|S_{q_k}|\geq 2$.\\

\noindent
(1) In this case, Lemma 7.2, (2) implies that if we apply ${\rm dav}_{fc}$ to $\mathbf{d}$ once, then every block of minimum value becomes empty. Namely we have ${\min }(\mathbf{d}^{(1)})\geq m+1$. Hence the inequality $w(\mathbf{d}^{(1)})<w(\mathbf{d})$ holds, which shows the validity of Proposition 6.2.\\

\noindent
(2) In this case, let $L\hspace{1mm} (\geq 1)$ denote the maximum value of $\ell_k, k\in [1,m]$. Then it follows from Lemma 7.2, (2) that if we apply ${\rm dav}_{fc}$ to $\mathbf{d}$ consecutively $L+1$ times, then every block of minimum value of $\mathbf{d}$ becomes empty. Namely we have ${\min }(\mathbf{d}^{(L+1)})\geq m+1$. Hence the inequality $w(\mathbf{d}^{(L)})<w(\mathbf{d})$ holds, which shows the validity of Proposition 6.2. \\

\noindent
(3) In this case, it follows from Lemma 7.2, (1) that if we apply ${\rm dav}_{fc}$ to $\mathbf{d}$ consecutively $\ell_k+1$ times, the length of the block $S^{(\ell_k+1)}_{q_k}$ of minimum value decreases by one from that of the original block $S_{q_k}$. Hence we have $\ell(\mathbf{d}^{(\ell_k+1)})<\ell(\mathbf{d})$, which shows the validity of Proposition 6.2. This completes the proof of Proposition 6.2. \qed\\

Combining Proposition 5.4 and Proposition 6.2, we complete the proof of the inductive step:

\begin{prp}
For any $\mathbf{d}\in\mathbf{Z}^n$ with $w(\mathbf{d})\geq 2$, the condition $C(\mathbf{d})$ holds.
\end{prp}

\section{Main theorem and its proof}
The following proposition lies almost near the goal:

\begin{prp}
For any $\mathbf{d}\in\mathbf{Z}^n$, there exists a nonnegative integer $K$ such that the width of $\mathbf{d}^{(K)}$ is equal to $0$ or $1$.
\end{prp}

\noindent
{\it Proof}.  We prove this by induction on the width. When $\mathbf{d}\in\mathbf{Z}^n$ itself has width 0 or 1, we can take $K=0$ and the proof finishes. Therefore we can assume that $w(\mathbf{d})=w_0\geq 2$ and that for any $\mathbf{d})$ with $w(\mathbf{d}')<w_0$ the assertion holds. We must show that there exists a positive integer $k$ such that 
\begin{eqnarray}
w(\mathbf{d}^{(k)})<w(\mathbf{d}).
\end{eqnarray}
Suppose, on the contrary, that the condition 
\begin{eqnarray}
w(\mathbf{d}^{(k)})=w(\mathbf{d})
\end{eqnarray}
for any positive integer $k$. Then it follows from Proposition * that there exists a positive integer $k_1$ such that 
\begin{eqnarray*}
\ell(\mathbf{d}^{(k_1)})<\ell(\mathbf{d})
\end{eqnarray*}
Furthermore applying Proposition * to $\mathbf{d}^{(k_1)}$, we see by () that there exists a positive integer $k_2$ such that
\begin{eqnarray*}
\ell(\mathbf{d}^{(k_1+k_2)})<\ell(\mathbf{d}^{(k_1)}).
\end{eqnarray*}
By repeating this argument, we obtain an increasing sequence $\{n_p\}$ of positive integers which makes an infinite decreasing sequence 
\begin{eqnarray*}
\ell(\mathbf{d}^{(n_1)})>\ell(\mathbf{d}^{(n_2)})>\ell(\mathbf{d}^{(n_3)})>\cdots,
\end{eqnarray*}
of positive integers, which is impossible. Thus the proof is completed. \qed\\

Combining Proposition 7.1 and Proposition 4.4 we finally complete the proof of our main theorem:

\begin{thm}
Let $Fix_D$ denote the set of the fixed points of the self map ${\rm dav}_{fc}:\mathbf{Z}^n\rightarrow\mathbf{Z}^n$, and let $Per_D$ denote the set of periodic points under ${\rm dav}_{fc}$. Then $Fix_D$ consists of the following two types of elements $\mathbf{d}=(d_0,\cdots,d_{n-1})$,
\begin{eqnarray*}
&&(1)\hspace{1mm}w(\mathbf{d})=0,\\
&&(2)\hspace{1mm}w(\mathbf{d})=1 \mbox{ and } d_i\in\{m,m+1\} \mbox{ for any }i\mbox{ with $m$ even},
\end{eqnarray*}
and $Per_D$ consists of the following elements,
\begin{eqnarray*}
&&(3)\hspace{1mm}w(\mathbf{d})=1 \mbox{ and } d_i\in\{m,m+1\} \mbox{ for any }i\mbox{ with $m$ odd}.
\end{eqnarray*}
The period of the elements of type $(3)$ is equal to $n$. Furthermore, for any $\mathbf{d}\in\mathbf{Z}^n\setminus (Fix_D\cup Per_D)$, there exists a positive integer $k$ such that ${\rm dav}_{fc}^k(\mathbf{d})\in Fix_D\cup Per_D$.
\end{thm}

We examine below what occurs if we apply $A$-average map to African good rhythms investigated in [1]. Here the parameters are fixed as $N=16, n=5$. It follows from Theorem 7.1 that $\mathbf{d}^{(k)}$ becomes periodic when its width equals one, namely when its entries consist solely of $3$ and $4$. For the presentation we introduce a distance to the cycle into which the orbit $\mathbf{d}^{(k)}\hspace{1mm},k\geq 0,$ falls:

\begin{df}
For any $\mathbf{d}\in {\rm CD}_N^{(n)}$, let $K$ be the minimum of the set $\{k\in\mathbf{Z}_{\geq 0};w(\mathbf{d}^{(k)})\leq 1\}$. We call $K$ the distance to the final cycle and denote it by ${\rm dist}_c(\mathbf{d})$.
\end{df}

\noindent
We present the African rhythms in the increasing order of the distance to the final cycle.\\

\noindent
Example 7.1.(1) ${\rm dist}_c(\mathbf{d})=0$:  "Bossa":$\mathbf{a}=(0,3,6,10,13)$. In this case $\mathbf{d}$ itself belongs to a cycle from the beginning: 
\begin{eqnarray*}
\begin{array}{lcl}
\mathbf{a}^{(0)}=(0,3,6,10,13) & \rightarrow & \mathbf{d}^{(0)}=(3,3,4,3,3) \\
\downarrow & & \downarrow\\
\vdots &  & \hspace{10mm}\vdots \\
\end{array}
\end{eqnarray*}

\noindent
(2.1) ${\rm dist}_c(\mathbf{d})=1$: "Shiko":$\mathbf{a}=(0,4,6,10,12)$. In this case $\mathbf{d}^{(k)}$ becomes periodic after $k=1$: 
\begin{eqnarray*}
\begin{array}{lcl}
\mathbf{a}^{(0)}=(0,4,6,10,12) & \rightarrow & \mathbf{d}^{(0)}=(4,2,4,2,4) \\
\downarrow & & \downarrow\\
\mathbf{a}^{(1)}=(2,5,8,11,14) & \rightarrow & \mathbf{d}^{(1)}=(3,3,3,3,4) \\
\downarrow & & \downarrow \\
\vdots &  & \hspace{10mm}\vdots \\
\end{array}
\end{eqnarray*}

\noindent
(2.2) ${\rm dist}_c(\mathbf{d})=1$: "Son":$\mathbf{a}=(0,3,6,10,12)$. In this case $\mathbf{d}^{(k)}$ becomes periodic after $k=1$: 
\begin{eqnarray*}
\begin{array}{lcl}
\mathbf{a}^{(0)}=(0,3,6,10,12) & \rightarrow & \mathbf{d}^{(0)}=(3,3,4,2,4) \\
\downarrow & & \downarrow\\
\mathbf{a}^{(1)}=(1,4,8,11,14) & \rightarrow & \mathbf{d}^{(1)}=(3,4,3,3,3) \\
\downarrow & & \downarrow \\
\vdots &  & \hspace{10mm}\vdots \\
\end{array}
\end{eqnarray*}

\noindent
(2.3) ${\rm dist}_c(\mathbf{d})=1$: "Rumba":$\mathbf{a}=(0,3,7,10,12)$. In this case $\mathbf{d}^{(k)}$ becomes periodic after $k=1$: 
\begin{eqnarray*}
\begin{array}{lcl}
\mathbf{a}^{(0)}=(0,3,7,10,12) & \rightarrow & \mathbf{d}^{(0)}=(3,4,3,2,4) \\
\downarrow & & \downarrow\\
\mathbf{a}^{(1)}=(1,5,8,11,14) & \rightarrow & \mathbf{d}^{(1)}=(4,3,3,3,3) \\
\downarrow & & \downarrow \\
\vdots &  & \hspace{10mm}\vdots \\
\end{array}
\end{eqnarray*}

\noindent
(3.1) ${\rm dist}_c(\mathbf{d})=3$: "Soukous":$\mathbf{a}=(0,3,6,10,11)$. In this case $\mathbf{d}^{(k)}$ becomes periodic after $k=3$: 
\begin{eqnarray*}
\begin{array}{lcl}
\mathbf{a}^{(0)}=(0,3,6,10,11) & \rightarrow & \mathbf{d}^{(0)}=(3,3,4,1,5) \\
\downarrow & & \downarrow\\
\mathbf{a}^{(1)}=(1,4,8,10,13) & \rightarrow & \mathbf{d}^{(1)}=(3,4,2,3,4) \\
\downarrow & & \downarrow \\
\mathbf{a}^{(2)}=(2,6,9,11,15) & \rightarrow & \mathbf{d}^{(2)}=(4,3,2,4,3) \\
\downarrow & & \downarrow \\
\mathbf{a}^{(3)}=(4,7,10,13,0) & \rightarrow & \mathbf{d}^{(3)}=(3,3,3,3,4) \\
\downarrow & & \downarrow \\
\vdots &  & \hspace{10mm}\vdots \\
\end{array}
\end{eqnarray*}

\noindent
(3.2) ${\rm dist}_c(\mathbf{d})=3$: "Gahu":$\mathbf{a}=(0,3,6,10,14)$. In this case $\mathbf{d}^{(k)}$ becomes periodic after $k=3$: 
\begin{eqnarray*}
\begin{array}{lcl}
\mathbf{a}^{(0)}=(0,3,6,10,14) & \rightarrow & \mathbf{d}^{(0)}=(3,3,4,4,2) \\
\downarrow & & \downarrow\\
\mathbf{a}^{(1)}=(1,4,8,12,15) & \rightarrow & \mathbf{d}^{(1)}=(3,4,4,3,2) \\
\downarrow & & \downarrow \\
\mathbf{a}^{(2)}=(2,6,10,13,0) & \rightarrow & \mathbf{d}^{(2)}=(4,4,3,3,2) \\
\downarrow & & \downarrow \\
\mathbf{a}^{(3)}=(4,8,11,14,1) & \rightarrow & \mathbf{d}^{(3)}=(4,3,3,3,3) \\
\downarrow & & \downarrow \\
\vdots &  & \hspace{10mm}\vdots \\
\end{array}
\end{eqnarray*}

\noindent
Note that the five rhythms (2.1)-(3.2) become rotationally equivalent to the first one "Bossa" after a number of iterations. Furthermore the rhythms with $w(\mathbf{d})=1$ fall into the category of "maximally even rhythms" investigated in [1].\\

\noindent
{\bf Reference}\\
\noindent
[1] Toussaint, G. {\it The Geometry of Musical Rhythm: What Makes a “Good” Rhythm Good?}, CRC Press, 2013.

\end{document}